\documentclass[12pt]{amsart}
\usepackage{geometry}
\usepackage{setspace}
\usepackage{enumerate,amssymb}
\usepackage{amsopn,amsfonts,amsmath,accents}
\usepackage{hyperref}
\hypersetup{colorlinks=true, pdfborderstyle={/S/U/W 1},
urlbordercolor={0 0 1}, linkcolor=red, urlcolor=blue,citecolor=gray}

\def\indic{{\rm {\large 1}\hspace{-2.3pt}{\large
l}}}
\usepackage{stackrel}
\allowdisplaybreaks

\def\In{\mbox{Int}}

\theoremstyle{definition}
\newtheorem{theorem}{Theorem}[section]
\newtheorem{corollary}{Corollary}[section]
\newtheorem{definition}{Definition}[section]

\newtheorem{ass}{Assumption}[section]
\numberwithin{equation}{section}
\allowdisplaybreaks

\begin{document}
\linespread{1.3}
\title[]{A triangular treatment effect model with random coefficients in the selection equation}
\author[Gautier]{Eric Gautier}
\address{Toulouse School of Economics, Universit\'e Toulouse 1 Capitole, 21 all\'ee de Brienne, 31000 Toulouse, France.}
\email{Eric.Gautier@tse-fr.eu}
\author[Hoderlein]{Stefan Hoderlein}
\address{Boston College, Chestnut Hill, MA 02467, USA.}
\email{Stefan.Hoderlein@bc.edu}
\date{This version: \today. 
This is a revision of 
\href{http://arxiv.org/abs/1109.0362v2}{arXiv:1109.0362v2}.} 
\thanks{We are grateful to seminar participants at Boston College, Bristol, Chicago, CREST, Harvard-MIT, Johns Hopkins, Kyoto, Nanterre, Northwestern, Oxford,Princeton, Toulouse, UCL, Vanderbilt, Yale, 2011 CIRM, 2012 Bates White, CLAPEM, ESEM, SCSE, World Congress in Probability and Statistics for useful comments, as well as Martin Browning, Kirill Evdokimov, James Heckman and Azeem Shaikh.}
\thanks{Eric Gautier acknowledges financial support from the grants ERC POEMH, ANR-13-BSH1-0004 and ANR-11-IDEX-0003/Labex Ecodec/ANR-11-LABX-0047}

\begin{abstract}
This paper considers treatment effects under endogeneity with complex heterogeneity in the selection equation. We model the outcome of an endogenous treatment as a triangular system, where both the outcome and first-stage equations consist of a random coefficients model. The first-stage specifically allows for nonmonotone selection into treatment. We provide conditions under which marginal distributions of potential outcomes, average and quantile treatment effects, all conditional on first-stage random coefficients, are identified. Under the same conditions, we derive bounds on the (conditional) joint distributions of potential outcomes and gains from treatment, and provide additional conditions for their point identification. All conditional quantities yield unconditional effects (\emph{e.g.}, the average treatment effect) by weighted integration.
\end{abstract}

\maketitle

\section{Introduction}\label{sec:0} 
To evaluate the effect of a treatment, it is common in social sciences to rely on non-experimental data.
In such a setup, ignoring self selection into treatment results in a misleading assessment of the effectiveness of the treatment, as causal effects may be confounded with the effect of endogenous selection.  Another important feature of such real-world applications is the heterogeneity of the effect of treatment across individuals.  The most prominent identification strategy relies on the (instrument) \emph{monotonicity} assumption of Imbens and Angrist (1994).  It means that there are \emph{no defiers} in the sense that
as an instrument is shifted from $z$ to $z'$ (say, a voucher to participate in treatment is randomly allocated), individuals either stay with their treatment choice or move into treatment, but no individual reverses his decision.  Vytlacil (2002) shows that this assumption is equivalent to modeling the selection as an additively separable latent index model with a single unobservable, where individuals are increasingly likely to be in the treatment group as the value of this unobservable increases. 

The monotonicity assumption is often too restrictive. The incentives 
individuals face might be complex and commonly cost-benefit trade-offs, as exemplified in the Roy model, are being performed. 
These considerations depend on several parameters, all of which may be expected to vary across the population. As such, more complex models of endogenous selection that explicitly allow for complex unobserved heterogeneity should be considered. A particularly natural specification for the selection equation that accomplishes this goal is given by a random coefficients
selection equation. Heckman and Vytlacil (2005) (HV henceforth) call it the benchmark nonseparable, nonmonotic model of treatment choice and emphasizes the importance of being equally flexible in terms of unobserved heterogeneity in the outcome and selection (or first-stage) equations.  In this paper, we consider a triangular system where the outcome equation is a linear random coefficients model with a binary endogenous regressor (\emph{i.e.}, the treatment status) and the selection equation is a nonlinear random coefficients binary choice model.  It allows for a population with both compliers and defiers.  Related literature includes  Beran and Millar (1994),  Beran, Feuerverger and Hall (1996) and Hoderlein, Klemel\"a and Mammen (2010) for the linear random coefficients model, 
Lewbel and Pendakur (2012) for a different nonlinear specification, Ichimura and Thomson (1998) and Gautier and Kitamura (2013) for the random coefficients binary choice model.
 
Under the monotonicity assumption, the Marginal Treatment Effect (MTE, Bj\"{o}rklund and Moffitt (1987) and HV (2005)) is a key structural parameter to recover a large variety of effects. 
In our framework without the monotonicity assumption, we show that marginal distributions of the potential outcomes, a generalization of the MTE and of the  Quantile Treatment Effect (QTE, Abadie, Angrist and Imbens (2002)), conditional on the first-stage random coefficients, are identified. 
Under the same assumptions, we obtain bounds on the joint distribution of potential outcomes and the distribution and variance of the gains, conditional on the first-stage random coefficients. We also provide conditions for their identification. The conditional effects that we obtain are structural parameters of interest in their own right, as we argue below. They are also building blocks for average, quantile or distributional 
effects. 

\section{The Theoretical Setup}\label{sec:1}
Throughout this paper, we make use of the following notations. The conditional distribution of a random vector $A$ given $B=b$ is denoted by $\mathbb{P}_{A|B}(\cdot|b)$; its conditional CDF by $F_{A|B}(\cdot|b)$; its conditional PDF by $f_{A|B}(\cdot|b)$ when it exists; its conditional expectation (resp. variance) by $\mathbb{E}[A|B=b]$ (resp. $Var(A|B=b)$) when $A$ is scalar; and the support of the conditional distribution by
$\mathrm{supp}\left(\mathbb{P}_{A|B}(\cdot|b)\right)$. The unconditional quantities are denoted without $|B$.
We use the notation $A\perp B|X$ when $A$ and $B$ are independent given $X$, a.s. for almost surely, a.e. for almost everywhere,
$\sigma (A)$ for the sigma algebra generated by $A$.  All random variables are defined on a probability space $(\Omega,\mathcal{F},\mathbb{P})$ and $\omega$ is an element of $\Omega$. 
We also denote by $\indic$ the indicator function, by $\In(A)$ the interior of a set $A$, by $\|\cdot\|$ the Euclidian norm, and by $A_{2-L}$ the subvector $(A_2,\dots,A_L)$  of $(A_1,\dots,A_L)$.

In the treatment effects framework, $Y_{0}$ is the outcome in the control group, $Y_{1}$ is the outcome in the treated group, $\Delta=Y_{1}-Y_{0}$ is the gain from treatment or treatment effect, $D$ is the binary variable which is equal to 1 when treatment has been selected  and 0 otherwise, and $Z$ is a vector of instruments.
We denote by $D_{z}(\omega)$ the value of $D(\omega)$ if  $Z(\omega)=(Z_1(\omega),\dots,Z_L(\omega))$ were externally set to $z$ (one has $D(\omega)=D_{Z(\omega)}(\omega)$). 

\subsection{Relaxing monotonicity}
Endogenous selection into treatment is often dealt with by specifying a selection equation with a single unobservable. In its general form, it can be written as
$D=\indic\left\{\mu(F(Z),\Theta)>0\right\}$,
where $Z$ is independent of $(Y_0,\Delta,\Theta)$, 
$F$ is measurable, $F(Z)$ and $\Theta$ are scalar, and $\mu$ is increasing in $F(Z)$. 
This equation can be written as an additively separable latent index model which, as shown in
Vytlacil (2002), is equivalent to the LATE assumptions of Imbens and Angrist (1994):
\begin{enumerate}[\textup{(}{L}1\textup{)}]
\item\label{L1} $\forall z\in\mathrm{supp}(\mathbb{P}_Z)$ $Z\perp (Y_0,Y_1,D_{z})$ and $z\to\mathbb{P}(D=1|Z=z)$ is not constant;
\item\label{L2} $\forall z,z'\in \mathrm{supp}(\mathbb{P}_Z)$,  $\forall \omega\in\Omega$ $D_{z}(\omega)\le D_{z'}(\omega)$ or $\forall \omega\in\Omega$ $D_{z'}(\omega)\le D_{z}(\omega)$.
\end{enumerate}
(L\ref{L2}) is called monotonicity.  
A few papers recently aimed at relaxing (L\ref{L2}). 
Small and Tan (2007) gives a formula for the Local Average Treatment Effect (LATE) as a function of quantities conditional on a latent variable $U$, but does not study their identification, and gives a degenerate example where $U=(Y_1,Y_0)$. Klein (2010) considers a selection equation with two unobservables and the order of the bias of classical estimators when one unobservable goes to zero.  Fox and Gandhi (2011) studies identification of the distribution of unobserved heterogeneity in Roy models with discretely supported random coefficients, without a random intercept in the selection equation. 

\subsection{The model} 
In this paper, we consider the triangular random coefficients model
\begin{align}
& Y=Y_0 +\Delta D,\label{m1} \\
& D=\indic\left\{ \Theta  -\Gamma_1(Z_1+g(Z_{2},\dots,Z_{L}))-\sum_{l=2}^{L}\Gamma_lf_l(Z_l)>0\right\}, \quad L\ge2  \label{m2}
\end{align}
where $f_2,\dots,\ f_{L}$ and $g$ are unknown. 
The random elements $Y$, $D$, $Z$ are observed, while $Y_0$, $\Delta$,  $\Theta$ and $\Gamma$ are unobserved.
The unobservables can be arbitrarily dependent. We do not assume the existence of densities or rely on parametric assumptions.

When available, conditioning variables $X$, with common values when $D$ is exogenously set to 0 or 1, can be introduced for two reasons: (1) to obtain effects for population subgroups (in which case they are not important in our identification arguments); (2) to justify the exogeneity of the instruments, or conditional independence assumptions. 

The outcome equation \eqref{m1} is a linear random coefficients equation with random coefficients $(Y_0,\Delta)$ and binary endogenous regressor $D$, namely $D$ and $(Y_0,\Delta)$ are dependent.
This specification imposes no restriction on the outcome equation.
The structural parameter in this model is $F_{Y_0,\Delta}$. It 
is in a one-to-one relationship with the joint distribution of potential outcomes $F_{Y_0,Y_1}$.
Important functionals of $F_{Y_0,\Delta}$ include: the average marginal effect $\mathbb{E}[\Delta]$, also called average treatment effect (ATE), the variance of treatment effects, and the proportion of people who benefit from treatment $\mathbb{P}(\Delta>0)$, among others.\\
\noindent{\bf Example RCR: A generalized Roy model with uncertainty and random coefficients.} 
Individuals select treatment if their \emph{ex-ante} gain from treatment $\mathbb{E}[\Delta|\mathcal{I}]$ exceeds cost.  The 
 \emph{ex-post} gain from treatment is $\Delta$. The sigma-field $\mathcal{I}$ is the information set individuals have at their disposal at the time of their decision to participate in the treatment.  Their 
cost is defined as $\Gamma_0+\Gamma_1(Z_1+g(Z_{2-L}))+\sum_{l=2}^{L}\Gamma_lf_l(Z_l)$. The instruments are cost shifters and the random coefficients reflect the heterogeneous impact of nonlinear transformations on cost.  We assume, for simplicity, that the individuals know the value of their random coefficients.  
We obtain model \eqref{m2} where $\Theta=\mathbb{E}[\Delta|\mathcal{I}]-\Gamma_0$.  The \emph{ex-ante} returns from treatment is  $\mbox{EAR}:=\mathbb{E}[\Delta|\mathcal{I}]-\Gamma_0-\Gamma_1(Z_1+g(Z_{2-L}))-\sum_{l=2}^{L}\Gamma_lf_l(Z_l)$, while the \emph{ex-post} returns is  $\mbox{EPR}:=\Delta-\Gamma_0-\Gamma_1(Z_1+g(Z_{2-L}))-\sum_{l=2}^{L}\Gamma_lf_l(Z_l)$\hfill$\square$

Model \eqref{m1}-\eqref{m2} nests the benchmark model put forward in HV (2005) and allows more general treatment choices than (L\ref{L2}) when $\Gamma_{2-L}$ has one nondegenerate coordinate. 

Because the individual values of a vector of random coefficients cannot be obtained when its dimension exceeds one, we cannot rely on a type of control function approach (see Newey, Powell and Vella (1999) and Imbens and Newey (2009)).  Kasy (2011) shows that, when the endogenous regressor is continuous, the control function approach cannot be used with a random coefficients first-stage. Identification of causal effect requires new tools, which are the  tools that we now develop. 

\subsection{Main Identifying Assumptions}\label{sec:11}
In this section, we present our main set of identifying assumptions.
Throughout this article, we make the following sign restriction.
\begin{ass}\label{assSign}
$Z\perp (\Theta,\Gamma)|X$ and $\Gamma_1$ has a sign and is not 0 a.s.
\end{ass}
Conditioning on $Z=z$ and $X=x$ for $(z,x)\in \In\left(\mathrm{supp}(\mathbb{P}_{Z,X})\right)$, increasing $z_1$ increases (resp. decreases) $\mathbb{P}(D=1|Z=z,X=x)$ if $\Gamma_1$ is negative (resp. positive). Therefore, the sign of $\Gamma_1$ can be identified. Since it is possible to change $Z_1$ in $-Z_1$, we now assume that $\Gamma_1>0$ a.s. We can then divide both sides of the inequality in \eqref{m2} by $\Gamma_1$ and work with the coefficients
$\overline{\Theta}:=\Theta/\Gamma_1$, $\overline{\Gamma}_2:=\Gamma_2/\Gamma_1,\dots$, $\overline{\Gamma}_L:=\Gamma_L/\Gamma_1$. 

The same argument shows that we do not have to know which instrument plays the role of $Z_1$. Indeed, it can be solely determined by $\mathbb{P}(D=1|Z=z,X=x)$. When $D$ is college attendance, tuition or distance to university could play the role of $Z_1$.

In the linear random coefficients binary choice model, assuming full support of the instruments and without invoking variables $X$, a sufficient condition for identification of the distribution of the random coefficients $(\Theta,\Gamma_1^0,\dots,\Gamma_L^0)$ up to scale is that, for a unit vector $s$, 
$s^T(\Theta,\Gamma_1^0,\dots,\Gamma_L^0)>0$ a.s. 
This means that that there exists a value $z$ of the instruments such that $D_z=1$ a.s. (in which case we take $s=(1,z)/\|(1,z)\|$) or $D_{z}=-1$ a.s. (in which case we take $s=-(1,z)/\|(1,z)\|$).  
This condition is satisfied if one coefficient has a sign. More generally, it means that by working with the vector of transformed instrument $AZ$ where $A$ is invertible and such that the first row of $(A^T)^{-1}$ is $s^T$, we obtain the new 
random coefficients $\Theta$ and $\Gamma=(A^T)^{-1}\Gamma^0$ where $\Gamma_1>0$ a.s. In equation \eqref{m2}, we replace the linear transformation by a nonlinear transformation, where $Z_1$ is transformed into $Z_1+g(Z_{2-L})$. 

We now invoke moments, instrument independence, and large support assumptions.
\begin{ass}\label{ass:Lu}
\begin{enumerate}[\textup{(}{A}1\textup{)}]
\item \label{A1} $\mathbb{E}[|Y_0|+|\Delta|]<\infty$;
\item \label{A2} $Z\perp (Y_d,\overline{\Theta},\overline{\Gamma}_{2-L})\ |X$ for $d=0,1$;
\item \label{A3} $0<\mathbb{P}(D=1|X)<1\quad a.s.$;
\item \label{A5}  $\forall(z_{2-L},x)\in \mathrm{supp}\left(\mathbb{P}_{Z_{2-L},X}\right)$, 
$$
\mathrm{supp}\left( \mathbb{P}_{Z_1|Z_{2-L},X}(\cdot |z_{2-L},x)\right) \supseteq
\mathrm{supp}\left( \mathbb{P}_{\overline{\Theta}-g(z_{2-L})-\sum_{l=2}^L\overline{\Gamma}_lf_l(z_l)|X}(\cdot |x)\right).$$
\end{enumerate}
\end{ass}
(A\ref{A1}) allows us to consider conditional expectations of $Y_0$ and $\Delta$.  (A\ref{A2}) allows $(Y_{0},\overline{\Theta},\overline{\Gamma}_{2-L})$ and $
(Y_{1},\overline{\Theta},\overline{\Gamma}_{2-L})$ to depend on $Z$, as long as we have at hand control variables $X$ which yield independence.  
These variables $X$ can either be observed directly, or they may be constructed from another additional equation, say, as control functions.  
Moreover, they need not exist at all.
(A\ref{A3}) states that there is a fraction of the
population that participates in treatment, and one that does not.
(A\ref{A5}) means that $Z_1$ has a large enough support. 
A similar assumption is made in Lewbel (2007). 
(A\ref{A5}) implies the existence of unselected samples when $Z_1$ (only) approaches the lower and upper bounds of its support, which could be tested.  
Traditional \textquotedblleft identification at infinity\textquotedblright makes inefficient use of the data, while  (A\ref{A5}) allows to build estimators that use the entire data (see Gautier and Hoderlein (2012)). Finally, though \textquotedblleft identification at infinity\textquotedblright arguments allow to obtain 
effects that depend on $F_{Y_d}$ for $d=0,1$, they
do not yield the structural parameters conditional on first-stage unobservables.

In a preliminary step, these assumptions may be taken to identify the nonlinear elements in \eqref{m2}. 
Observe that using (A\ref{A2}), for every $(z,x)\in \mathrm{supp}\left(\mathbb{P}_{Z,X}\right)$,
$$\mathbb{P}\left(D=0|Z=z,X=x\right)=\mathbb{P}\left(
\overline{\Theta}-g(z_{2-L})-\sum_{l=2}^L\overline{\Gamma}_lf_l(z_l)\le z_1|X=x\right)\ .$$
Therefore, using (A\ref{A5}), identification in \eqref{m2} is implied by identification in
$$W=g(Z_{2-L})+\sum_{l=2}^L\overline{\Gamma}_lf_l(Z_l)-\overline{\Theta}$$
where $Z_{2-L}\perp (\overline{\Gamma}_{2-L},\overline{\Theta}) |X$. We now make integrability and location normalizations.
\begin{ass}\label{ass:nonlinear}
\begin{enumerate}[\textup{(}{N}1\textup{)}]
\item \label{N0}  $\mathbb{E}\left[|g(Z_{2-L})|+\sum_{l=2}^L|\overline{\Gamma}_{l}|+|f_l(Z_l)|\right]<\infty$;
\item \label{N1}  For $l=2,\dots,L$, $\mathbb{E}[g(Z_{2-L})|Z_l,X]=0$ a.s.;
\item \label{N2} For $l=2,\dots,L$, $\mathbb{E}[\overline{\Gamma}_l|X]=1$ a.s.;
\item \label{N3} For $l=2,\dots,L$, $\mathbb{E}\left[f_l(Z_l)|X\right]=0$ a.s. 
\end{enumerate}
\end{ass}
The following theorem shows identification of the nonlinear elements in \eqref{m2}.
\begin{theorem}\label{thrm:functions}
$f_l$ for $l=2,\dots,L$ and $g$ are identified under Assumption \ref{ass:nonlinear}.
\end{theorem}
Writing $\overline{Z}_1:=-Z_1+g(Z_{2-L})$ and $\overline{Z}_l:=f_l(Z_l)$ for $l=2,\dots,L$, we now work with
\begin{equation}
D=\indic\left\{ -\overline{\Theta}+\sum_{l=2}^{L}\overline{\Gamma}_l\overline{Z}_l<\overline{Z}_1\right\} \label{m2b}
\end{equation}
where $\mathbb{E}[\overline{\Gamma}_l|X]=1$ for $l=2,\dots,L$.

When the large support assumption (A\ref{A5}) is not satisfied, it is replaced by (A4'):

\begin{enumerate}[\textup{(}{A4'}a\textup{)}]
\item Model \eqref{m2} does not have unknown functions;
\item \label{SLS1} $\forall x\in\mathrm{supp}\left(\mathbb{P}_X\right)$, $\forall R\ge 0$, $
\mathbb{E}\left[\left.\exp\left(R|\overline{\Theta}|\right)\right| X=x\right]<\infty$;
\item \label{SLS2} $\forall (z_{2-L},x)\in\mathrm{supp}\left(\mathbb{P}_{\overline{Z}_{2-L},X}\right)$, $\In\left(\mathrm{supp}\left(\mathbb{P}_{Z_1|\overline{Z}_{2-L},X}(\cdot|z_{2-L},x)\right)\right)$ is not empty.
\end{enumerate}

Finally, we invoke either of the two following assumptions on ${ \rm supp}(\mathbb{P}_{\overline{Z}_{2-L}|X})$.
\begin{ass}\label{ass:I}
$\forall x\in \mathrm{supp}\left(\mathbb{P}_X\right)$,
$\mathrm{supp}\left(\mathbb{P}_{\overline{Z}_{2-L}|X}(\cdot|x)\right)=\mathbb{R}^{L-1}$.
\end{ass}

\begin{ass}\label{ass:Ib}
$\forall x\in\mathrm{supp}\left(\mathbb{P}_X\right)$,
\begin{enumerate}[\textup{(}{S}1\textup{)}]
\item \label{S1} $\forall R\ge0$, 
$\mathbb{E}\left[\left.\exp\left(R\left\|\overline{\Gamma}_{2-L}\right\|\right)\right| X=x\right]<\infty$;
\item \label{S2} $\In\left(\mathrm{supp}\left(\mathbb{P}_{\overline{Z}_{2-L}|X}(\cdot|x)\right)\right)$ is not empty.
\end{enumerate}
\end{ass}
(S\ref{S1}) corresponds to the existence of a moment generating function of $\overline{\Gamma}_{2-L}$, which can be relaxed as in Gaillac and Gautier (2015). Beran and Millar (1994) assumes 
that the support of the random coefficients, including the random intercept, is compact.

\section{Parameters Based on the Marginals of Potential Outcomes}\label{sec:2}
\subsection{Identification of Marginals Conditional on First-stage Unobservables}\label{sec:21}
The following result is central in our analysis.
\begin{theorem}\label{thrm:1}
Under Assumption \ref{ass:Lu} with (A\ref{A5}) or (A4') and either Assumption \ref{ass:I} or \ref{ass:Ib}, $F_{Y_{d},\overline{\Theta},\overline{\Gamma}_{2-L},X}$ for $d=0,1$ are identified.
\end{theorem}
This result shows identification of $F_{Y_0+\Delta d|\overline{\Theta},\overline{\Gamma}_{2-L},X}$ and $F_{Y_0+\Delta d|X}$ for $d=0,1$, and of the average and quantile structural functions.
This means that we have resolved the endogeneity issue in model \eqref{m1} entirely and the difficulty in identifying the structural parameter $F_{Y_0,\Delta|X}$ now stems solely from the fact that the regressor is binary. 

As will be obvious, Theorem \ref{thrm:1} may be used to establish identification of many parameters. For instance, it directly yields identification of the following parameter.
\begin{definition} The Unobservables Conditioned Quantile Treatment Effect (UCQTE) is defined, for $(\theta,\gamma_{2-L},x)$ in $\mathrm{supp}(\mathbb{P}_{\overline{\Theta},\overline{\Gamma}_{2-L},X})$, as
\begin{equation*}
\mathrm{UCQTE}(\theta,\gamma_{2-L},x,\tau):=F_{Y_{1}|\overline{\Theta},\overline{\Gamma}_{2-L},X}(\cdot|\theta,\gamma_{2-L},x)^{-1}(\tau)-F_{Y_{0}|\overline{\Theta},\overline{\Gamma}_{2-L},X}(\cdot|\theta,\gamma_{2-L},x)^{-1}(\tau)\ .
\end{equation*}
\end{definition}
Integrating out the first-stage unobserved heterogeneity, we obtain 
\begin{equation*}
{\rm QTE}(x,\tau )=F_{Y_{1}|X}(\cdot|x)^{-1}(\tau)-F_{Y_{0}|X}(\cdot|x)^{-1}(\tau),\ \tau\in(0,1)
\end{equation*}
without monotonicity or rank invariance (Chernozhukov and Hansen (2005)).  

\subsection{The Unobservables Conditioned Average Treatment Effect (UCATE)}\label{sec:UCATE}
\subsubsection{Definition and Properties}
MTE is a key structural parameter 
(see HV (2005)) which relies on the additively separable latent index model for the first-stage, thus on monotonicity.
To extend this concept to our framework, 
we introduce the following generalization of the MTE. 
It is clearly identified from Theorem \ref{thrm:1}.
\begin{definition}
The Unobservables Conditioned Average Treatment Effect (UCATE) is defined, for $(\theta,\gamma_{2-L},x)$ in $\mathrm{supp}(\mathbb{P}_{\overline{\Theta},\overline{\Gamma}_{2-L},X})$, as
\begin{equation*}
\mathrm{UCATE}(\theta,\gamma_{2-L},x):=\mathbb{E}[\Delta |\overline{\Theta}=\theta,\overline{\Gamma}_{2-L}=\gamma_{2-L},X=x]\ .
\end{equation*}
\end{definition}

UCATE shares the same interpretation and attractive properties as MTE. It is the average effect for a subpopulation with $X=x$ which would be indifferent between participation and
nonparticipation in the treatment, if it were exogenously
assigned a value $z$ of $\overline{Z}$ such
that $-\overline{\Theta}+\sum_{l=2}^{L}\overline{\Gamma}_lz_l=z_1$.  Due to (A\ref{A2}), this parameter is policy invariant (\emph{i.e.}, it is independent of the values of the instruments)
$$\mathbb{E}\left[\left. \Delta \right|\overline{\Theta},\overline{\Gamma}_{2-L},X,\overline{Z}\right]=\mathbb{E}[\left.\Delta \right|\overline{\Theta},\overline{\Gamma}_{2-L},X]\ .$$
Like MTE, UCATE is an economically important parameter. 
To fix ideas, think of the treatment as attending college,
and of the $j$-th cost factor as college tuition. If an individual has a
high value of $\overline{\Gamma}_{l}$,
she is responsive to a change in college tuition. This may be because she is liquidity constrained.  For a policy maker who is deciding about tuition rates, 
the difference in treatment effects across these groups is a key variable in assessing the effect of a change in the incentive structure. For targeted policy measures, it is interesting to see how the effect of treatment varies with unobserved sensitivity to tuition and 
observed $X$ ({\em e.g.}, race). 

Finally, note that, if $\mathrm{UCATE}(\theta,\gamma_{2-L},x)$ is not constant in one of the $\gamma_l$'s, it is an indication that monotonicity does not hold. This allows to test monotonicity.

\subsubsection{Example RCR: UCATE as the Ex-ante Gains From Treatment}\label{sec:UCATEexE}
In this model, UCATE   
is a core structural effect. We invoke the following assumption.
\begin{ass}\label{ass:LDh}
 \begin{enumerate}[\textup{(}{U}1\textup{)}]
\item \label{U1} $\mathcal{I}\supseteq \sigma(\overline{\Gamma}_{2-L},X)$
\item \label{U2} $\sigma(\overline{\Gamma}_{2-L},X)=\sigma(\Gamma_0,\Gamma_1,\dots,\Gamma_{L},X)$\ .
\end{enumerate}
\end{ass}
Under Assumption \ref{ass:LDh}, individuals have perfect knowledge of their cost function, but can be uncertain about the potential outcomes. 

(U\ref{U2}) is satisfied when $(\Gamma_0,\Gamma_1,\dots,\Gamma_{L})=\varphi(\overline{\Gamma}_{2-L})$ for some measurable function $\varphi$.  This restricts the unobserved heterogeneity entering \eqref{m2}. However, it does not restrict its dimension. Example RCR also involves $\mathbb{E}[\Delta|\mathcal{I}]$ in \eqref{m2}, so that, when the superset symbol in (U\ref{U1}) is strict, \eqref{m2} depends on a nondegenerate vector of unobservables of dimension $L$.  A sufficient condition for (U\ref{U2}) to hold is that one $\Gamma_l$ for $l=1,\dots,L$ is $\sigma(X)$-measurable and $\mathbb{P}(\Gamma_l=0)=0$, and $\Gamma_0$ is $\sigma(\Gamma_1,\dots,\Gamma_L,X)$-measurable. Still, it is a richer specification than the extended Roy model where the cost function is not random.

Define  the errors $\eta_d$ for $d=0,1$ as
\begin{equation}\label{eq:info1}
\eta_d:=\mathbb{E}[Y_d|\mathcal{I}]-\mathbb{E}[Y_d|\overline{\Gamma}_{2-L},X]\ .
\end{equation}
When Example RCR models sectorial choice, they could be interpreted as sector specific skills known by the individuals at the time the choice is made and unobserved by the econometrician.
Because $\overline{\Theta}=(\mathbb{E}[Y_1-Y_0|\mathcal{I}]-\Gamma_0)/\Gamma_1=(\mathbb{E}[Y_1-Y_0|\overline{\Gamma}_{2-L},X]+\eta_1-\eta_0-\Gamma_0)/\Gamma_1$,  Assumption \ref{ass:LDh} implies
\begin{equation}\label{eq:info3}
\sigma(\overline{\Theta},\overline{\Gamma}_{2-L},X)=
\sigma(\overline{\Gamma}_{2-L},\eta_1-\eta_0,X)\subseteq \sigma(\overline{\Gamma}_{2-L},\eta_1,\eta_0,X)\subseteq \mathcal{I}\ .
\end{equation}
Even if $\sigma(\overline{\Theta},\overline{\Gamma}_{2-L},X)$ is a strict subset of $\mathcal{I}$, the following result holds.
\begin{theorem}\label{prop:ucate}
Under Assumption \ref{ass:LDh} with (A\ref{A5}) or  (A4'), $\mathrm{UCATE}(\overline{\Theta},\overline{\Gamma}_{2-L},X)=\mathbb{E}[\Delta|\mathcal{I}]$. 
Moreover, if $\mathbb{E}[\Delta|\mathcal{I}]$ is not $\sigma(\overline{\Gamma}_{2-L},X)$-measurable, $F_{\left.\mbox{EAR}\right|Z,X}$ is identified.
\end{theorem}

\subsubsection{UCATE as a Building Block to Obtain Treatment Effects that Depend on Averages}
Like MTE,  UCATE is a building block to obtain 
treatment effect parameters that depend on averages. For example, the ATE and Treatment
on the Treated (TT) are given by
\begin{align*}
\mathrm{ATE}(x)& =\int_{\mathbb{R}^{L}}\mathrm{UCATE}(\theta,\gamma_{2-L},x)
d\mathbb{P}_{\overline{\Theta},\overline{\Gamma}_{2-L}|X}(\theta,\gamma_{2-L}|x)\ ,\\
\mathrm{TT}(x)& =\int_{\mathbb{R}^{L}}h_{\mathrm{TT}}(\theta,\gamma_{2-L},x)
\mathrm{UCATE}(\theta,\gamma_{2-L},x)d\mathbb{P}_{\overline{\Theta},\overline{\Gamma}_{2-L}|X}(\theta,\gamma_{2-L}|x)\ ,
\end{align*}
where $h_{\mathrm{TT}}(\theta,\gamma_{2-L} ,x)=\mathbb{E}\left[ D |X=x\right]
^{-1}\mathbb{E}\left[ \left.
\indic\left\{ -\theta+\sum_{l=2}^{L}\gamma_l\overline{Z}_l<\overline{Z}_1\right\}
\right |X=x\right] $.

UCATE also allows to derive Policy Relevant Treatment Effect (PRTE) parameters. In our setup where an individual self-select into treatment, these parameters inform the decision maker about the effect of a change in
the structure of the variables or incentives she can control, namely the distribution of instruments. For example, the alternative policy can
consist in changing $Z_1$ to $Z_1+\alpha$ (\textit{e.g.}, increasing tuition by $\alpha$).
HV (2001) formalizes this idea. 
Denote by $\mathbb{E}_{\alpha}$ the expectation under
a policy $\alpha$ and $\mathbb{E}_{0}$ the expectation under the baseline policy. 
The PRTE parameter is defined, for $x$ in $\mathrm{supp}(\mathbb{P}_{X})$, as
$$\mathrm{PRTE}(\alpha,x):=\frac{\mathbb{E}_{\alpha}[Y|X=x]-\mathbb{E}_{0}[Y|X=x]}{\mathbb{E}_{\alpha}[D|X=x]-\mathbb{E}_{0}[D|X=x]}\ .$$
Marginal Policy Relevant Treatment Effects (MPRTE, Carneiro, Heckman and Vytlacil (2010))
are defined as 
$$
\mathrm{MPRTE}(x):=\lim_{\alpha\to 0}\mathrm{PRTE}(\alpha,x)\ .
$$
In our setup,  
$\mathrm{PRTE}(\alpha,x)$ can be obtained from UCATE like $\mathrm{TT}(x)$ replacing the weight $h_{\mathrm{TT}}(\theta,\gamma_{2-L},x)$ by 
$h_{\mathrm{PRTE}}(\theta,\gamma_{2-L},\alpha,x)$ defined by
$$\frac{\mathbb{E}_{\alpha}\left[ \left.
\indic\left\{ -\theta+\sum_{l=2}^{L}\gamma_l\overline{Z}_l<\overline{Z}_1\right\}
\right |X=x\right]-\mathbb{E}_0\left[ \left.
\indic\left\{ -\theta+\sum_{l=2}^{L}\gamma_l\overline{Z}_l<\overline{Z}_1\right\}
\right |X=x\right]}{\mathbb{E}_{\alpha}\left[\left.
\indic\left\{ -\overline{\Theta}+\sum_{l=2}^{L}\overline{\Gamma}_l\overline{Z}_l<\overline{Z}_1\right\}
\right |X=x\right]-\mathbb{E}_{0}\left[ \left.
\indic\left\{ -\overline{\Theta}+\sum_{l=2}^{L}\overline{\Gamma}_l\overline{Z}_l<\overline{Z}_1\right\}
\right |X=x\right]}\ .
$$
In a similar manner, $\mathrm{MPRTE}(x)$ can be obtained using the weight
$$h_{\mathrm{MPRTE}}(\theta,\gamma_{2-L},x)=\frac{\partial_{\alpha}\left.\mathbb{E}_{\alpha}\left[ \left.
\indic\left\{ -\theta+\sum_{l=2}^{L}\gamma_l\overline{Z}_l<\overline{Z}_1\right\}
\right |X=x\right]\right|_{\alpha=0}}{\partial_{\alpha}\left.\mathbb{E}_{\alpha}\left[ \left.
\indic\left\{ -\overline{\Theta}+\sum_{l=2}^{L}\overline{\Gamma}_l\overline{Z}_l<\overline{Z}_1\right\}
\right |X=x\right]\right|_{\alpha=0}}\ .
$$
The weights $h_{TT}$, $h_{PRTE}$ and $h_{MPRTE}$ can be obtained by Monte-Carlo techniques. 

\section{Parameters Based on the Joint Distribution of Potential Outcomes}\label{sec:3}
In this section, we consider identification of the structural parameter $F_{Y_0,\Delta|X}$ in model \eqref{m2}.  The importance of distributional treatment effects is emphasized by Heckman, Smith and Clements (1997) (HSC henceforth) and Abbring and Heckman (2007). We present partial identification results under the same assumptions as in Section \ref{sec:2} and present assumptions that allow to identify distributional effects. 

\subsection{Partial Identification of Distributional Effects}\label{sec:32}
Theorem \ref{thrm:1} yields that, under our maintained assumptions, 
$F_{Y_d|\overline{\Theta},\overline{\Gamma}_{2-L},X}$ for $d=0,1$
are identified. Under the same assumptions, we obtain sharp bounds on structural parameters that depend on 
$F_{Y_0,Y_1|\overline{\Theta},\overline{\Gamma}_{2-L},X}$. We use the following notations:
\begin{align*}
F_{Y_0,Y_1|\overline{\Theta},\overline{\Gamma}_{2-L},X}^{L}(y_{0},y_{1}|\theta,\gamma_{2-L},x)
& :=\max \left\{F_{Y_0|\overline{\Theta},\overline{\Gamma}_{2-L},X}(y_0|\theta,\gamma_{2-L},x)+F_{Y_1|\overline{\Theta},\overline{\Gamma}_{2-L},X}(y_1|\theta,\gamma_{2-L},x)-1,0\right\}\\
F_{Y_0,Y_1|\overline{\Theta},\overline{\Gamma}_{2-L},X}^{U}(y_{0},y_{1}|\theta,\gamma_{2-L},x)
& :=\min \left\{F_{Y_0|\overline{\Theta},\overline{\Gamma}_{2-L},X}(y_0|\theta,\gamma_{2-L},x),F_{Y_1|\overline{\Theta},\overline{\Gamma}_{2-L},X}(y_1|\theta,\gamma_{2-L},x)\right\}
\end{align*}
$$F_{\Delta|\overline{\Theta},\overline{\Gamma}_{2-L},X}^{L}(\delta |\theta,\gamma_{2-L},x) :=\sup_{y\in \mathbb{R}}\max \left\{
F_{Y_1|\overline{\Theta},\overline{\Gamma}_{2-L},X}(y|\theta,\gamma_{2-L},x)-F_{Y_0|\overline{\Theta},\overline{\Gamma}_{2-L},X}(y-\delta|\theta,\gamma_{2-L},x),0\right\}$$
$$F_{\Delta|\overline{\Theta},\overline{\Gamma}_{2-L},X}^{U}(\delta |\theta,\gamma_{2-L},x) :=1+\inf_{y\in \mathbb{R}}\min \left\{
F_{Y_1|\overline{\Theta},\overline{\Gamma}_{2-L},X}(y|\theta,\gamma_{2-L},x)-F_{Y_0|\overline{\Theta},\overline{\Gamma}_{2-L},X}(y-\delta|\theta,\gamma_{2-L},x),0\right\}
$$
and $F_{Y_0,Y_1}^{L}(y_{0},y_{1}|x)$, $F_{Y_0,Y_1}^{U}(y_{0},y_{1}|x)$, $F_{\Delta}^{L}(\delta|x)$, and $F_{\Delta}^{U}(\delta|x)$ are the expectation of above quantities evaluated at the random vectors. For example, we define 
$$
F_{Y_0,Y_1}^{L}(y_{0},y_{1}|x) :=\int_{\mathbb{R}^{L}}F_{Y_0,Y_1|\overline{\Theta},\overline{\Gamma}_{2-L},X}^{L}(y_{0},y_{1}|\theta,\gamma_{2-L},x)
d\mathbb{P}_{\overline{\Theta},\overline{\Gamma}_{2-L}|X}(\theta,\gamma_{2-L}|x)\ .$$

\begin{theorem}\label{theorem:BoundsMU} Under the assumptions of Theorem \ref{thrm:1}, for $(\theta,\gamma_{2-L},x)\in\mathrm{supp}(\mathbb{P}_{\overline{\Theta},\overline{\Gamma}_{2-L},X})$, $(y_{0},y_{1},\delta )\in \mathbb{R}^{3}$,
\begin{equation*}
F_{Y_0,Y_1|\overline{\Theta},\overline{\Gamma}_{2-L},X}^{L}(y_{0},y_{1}|\theta,\gamma_{2-L},x)\leq F_{Y_{0},Y_{1}|\overline{\Theta},\overline{\Gamma}_{2-L},X}(y_{0},y_{1}|\theta,\gamma_{2-L},x)\leq
F_{Y_0,Y_1|\overline{\Theta},\overline{\Gamma}_{2-L},X}^{U}(y_{0},y_{1}|\theta,\gamma_{2-L},x)
\end{equation*}
\begin{equation*}
F_{\Delta|\overline{\Theta},\overline{\Gamma}_{2-L},X}^{L}(\delta |\theta,\gamma_{2-L},x)\leq F_{\Delta|\overline{\Theta},\overline{\Gamma}_{1-2},X}(\delta |\theta,\gamma_{2-L},x)\leq
F_{\Delta|\overline{\Theta},\overline{\Gamma}_{2-L},X}^{U}(\delta |\theta,\gamma_{2-L},x)\ .
\end{equation*}
\end{theorem}

All bounds are sharp. Bounds on functionals of the unobservables conditioned joint distribution of potential outcomes can be obtained like in Fan, Guerre and Zhu (2014).

\begin{corollary}\label{cor:1} Under the assumptions of Theorem \ref{thrm:1}, for $x\in\mathrm{supp}(\mathbb{P}_{X})$, $(y_{0},y_{1},\delta )\in \mathbb{R}^{3}$,
\begin{equation*}
F_{Y_0,Y_1|X}^{L}(y_{0},y_{1}|x)\leq F_{Y_{0},Y_{1}|X}(y_{0},y_{1}|x)\leq
F_{Y_0,Y_1|X}^{U}(y_{0},y_{1}|x)
\end{equation*}
\begin{equation*}
F_{\Delta|X}^{L}(\delta |x)\leq F_{\Delta|X}(\delta |x)\leq
F_{\Delta|X}^{U}(\delta |x)\ .
\end{equation*}
\end{corollary}
Bounds on the above quantities have been obtained in the case of randomized experiments or selection on observables (see HSC (1997) for $F_{Y_0,Y_1}$, Fan and  Park (2010) and Firpo and Ridder (2008) for $F_{\Delta}$).
The bounds of this section hold when there is endogenous selection into treatment and without mononoticity. 
It is possible to deduce bounds which are unconditional on (some) $X$ by integration against $\mathbb{P}_X$. These are sharper than the ones obtained without conditioning (see Firpo and Ridder (2008)). 
Similarly, the bounds of Corollary \ref{cor:1} are sharper than bounds obtained without monotonicity (see Kitagawa (2009) for $F_{Y_0,Y_1}$) and without conditioning on the first-stage unobservables. 

\subsection{Point Identification of Distributional Effects}\label{sec:33}
We start by strengthening 
(A\ref{A2}):
\begin{center}
(A2')\hspace{0.5cm} $Z\perp (Y_0,Y_1,\overline{\Theta},\overline{\Gamma}_{2-L})\ |X$. 
\end{center}
We now present three alternative identifying assumptions and give sufficient conditions for them to hold in the case of Example RCR under Assumption \ref{ass:LDh}. 

\begin{ass}\label{ass:condindep1}
$Y_0 \perp Y_1\ |\overline{\Gamma}_{2-L},X$ .
\end{ass}
\begin{ass}\label{ass:condindep2}
$Y_0 \perp Y_1\ |\overline{\Theta},\overline{\Gamma}_{2-L},X$ .
\end{ass}
\begin{ass}\label{ass:condindep3}
\begin{enumerate}[\textup{(}i\textup{)}]
\item\label{assi} $Y_0\perp\Delta\ |\overline{\Theta},\overline{\Gamma}_{2-L},X$; 
\item\label{assii} Either (1) for a.e. $t\in\mathbb{R}$, $\mathbb{E}[e^{itY_0}| \overline{\Theta},\overline{\Gamma}_{2-L},X]\ne 0$ a.s.\\ or (2) $\forall R\ge 0$, 
$
\mathbb{E}\left[\left.
\exp\left(R
|\Delta|\right)
\right| \overline{\Theta},\overline{\Gamma}_{2-L},X\right]<\infty
$ a.s.
\end{enumerate}
\end{ass}

HSC (1997) makes a similar assumption as Assumption \ref{ass:condindep3} \eqref{assi}, but it is not conditional on the first-stage unobservables. 
Condition \eqref{assii} (1) is the condition in Devroye (1989) for the deconvolution problem $Y_1=\Delta+Y_0$ where $Y_0$ is the error with known distribution. Condition \eqref{assii} (2) imposes no restriction on the distribution of $Y_0$.

We can relabel state 0 as state 1 in the first two assumptions, as well as in Assumption  \ref{ass:condindep3} \eqref{assi} if the agents have perfect foresight of the outcomes in both states.  This is important when the two states are two sectors in the economy.  

Consider now Example RCR. Define $\epsilon_d:=Y_d-\mathbb{E}[Y_d|\mathcal{I}]$ for $d=0,1$, which could be interpreted as state specific unexpected shocks. Consider the three following assumptions:
\begin{enumerate}[\textup{(}{B}1\textup{)}]
\item\label{identdistrib1} $\eta_0+\epsilon_0 \perp \eta_1+\epsilon_1\ |\overline{\Gamma}_{2-L},X$;
\item\label{identdistrib2} Assumption \ref{ass:LDh} holds, $\epsilon_0 \perp \epsilon_1\ |\overline{\Gamma}_{2-L},\eta_1-\eta_0,X$, and
$\eta_0$ is  $\sigma(\overline{\Gamma}_{2-L},\eta_1-\eta_0,X)$ measurable;
\item\label{identdistrib3}  Assumption \ref{ass:LDh} holds and $\eta_0+\epsilon_0 \perp \epsilon_1-\epsilon_0\ |\overline{\Gamma}_{2-L},\eta_1-\eta_0,X$.
\end{enumerate}

(B\ref{identdistrib1}) and (B\ref{identdistrib2}) 
restrict the unobserved heterogeneity in \eqref{m2}. (B\ref{identdistrib1}) implies
$$Y_0-\mathbb{E}[Y_0 | \overline{\Gamma}_{2-L},X] \perp Y_1-\mathbb{E}[Y_1 | \overline{\Gamma}_{2-L},X]\ |\overline{\Gamma}_{2-L},X\ ,$$
and thus Assumption \ref{ass:condindep1}.
(B\ref{identdistrib2}) implies
$$Y_0-\mathbb{E}[Y_0 | \overline{\Theta},\overline{\Gamma}_{2-L},X] \perp Y_1-\mathbb{E}[Y_1 | \overline{\Theta},\overline{\Gamma}_{2-L},X]\ |\overline{\Theta},\overline{\Gamma}_{2-L},X$$
and thus Assumption \ref{ass:condindep2}. (B\ref{identdistrib3}) implies Assumption \ref{ass:condindep3} \eqref{assi}. It is the less restrictive on the unobserved heterogeneity entering the first-stage equation. It is satisfied, for example, when there is no uncertainty on the outcome in the base state (\emph{i.e.}, $\epsilon_0=0$) and $\epsilon_1\perp \eta_0, \eta_1,\overline{\Gamma}_{2-L}|X$. 
Assumptions \ref{ass:condindep1} - \ref{ass:condindep3} \eqref{assi} can hold more generally, in particular without Assumption \ref{ass:LDh} which restricts the unobserved heterogeneity entering \eqref{m2}.

The next theorem shows point identification under either of the three assumptions.
\begin{theorem}\label{thrm:CondIndep}
Under the assumptions of Theorem \ref{thrm:1} replacing (A\ref{A2}) by (A2'), $F_{Y_0,Y_1|\overline{\Gamma}_{2-L},X}$ is identified under
Assumption \ref{ass:condindep1}, while $F_{Y_0,Y_1|\overline{\Theta},\overline{\Gamma}_{2-L},X}$ is identified under Assumption \ref{ass:condindep2} or \ref{ass:condindep3}. 
\end{theorem}

The following result is consequence Theorem \ref{prop:ucate}.
 Carneiro, Hansen and Heckman (2003) and Cunha and Heckman (2007) rely on a factor model to obtain related results.
\begin{theorem}\label{thrm:ex} Suppose \eqref{m2} corresponds to model RCR. 
Under the assumptions of Theorem \ref{thrm:CondIndep} with Assumption \ref{ass:condindep2} or \ref{ass:condindep3}, and Assumption \ref{ass:LDh},
$F_{\Delta,\mathbb{E}[\Delta|\mathcal{I}]|X}$ is identified.
Moreover, if $\mathbb{E}[\Delta|\mathcal{I}]$ is not $\sigma(\overline{\Gamma}_{2-L},X)$-measurable, then 
$F_{\left.\mbox{EAR},\mbox{EPR}\right|Z,X}$
is identified.
\end{theorem}

UCATE is a building block to obtain effects that depend on averages. The structural parameter that we now introduce plays a similar role, but allows to obtain all effects that depend on $F_{Y_0,Y_1|X}$. 
It is clearly identified from Theorem \ref{thrm:CondIndep} and policy invariant.
\begin{definition} The Unobservables Conditioned Distribution of Treatment Effects is defined, for $\delta$ in $\mathbb{R}$ and $(\theta,\gamma_{2-L},x)$ in $\mathrm{supp}(\mathbb{P}_{\overline{\Theta},\overline{\Gamma}_{2-L},X})$, as
\begin{equation*}
\mathrm{UCDITE}(\delta ,\theta,\gamma_{2-L} ,x):=f_{\Delta |\overline{\Theta},\overline{\Gamma}_{2-L},X}(\delta
|\theta,\gamma_{2-L},x)\ .
\end{equation*}
\end{definition}
UCDITE is a key quantity to obtain all effects that  depend on $F_{Y_0,Y_1|X}$, 
for example
\begin{align*}
\mathbb{P}\left( \Delta >0|X=x\right) & =\int_{\mathbb{R}}\indic\left\{
\delta >0\right\} \int_{\mathbb{R}^{L}}{\rm UCDITE}(\delta ,\theta,\gamma_{2-L} ,x)
d\mathbb{P}_{\overline{\Theta},\overline{\Gamma}_{2-L}|X}(\theta,\gamma_{2-L}|x)\\
f_{Y_{0},Y_{1}|X}(y_{0},y_{1}|x)& =\int_{\mathbb{R}^{L}}
{\rm UCDITE}(y_{1}-y_{0},\theta,\gamma_{2-L} ,x)d\mathbb{P}_{Y_0,\overline{\Theta},\overline{\Gamma}_{2-L}|X}(y_0,\theta,\gamma_{2-L}|x) \\
f_{\Delta |D=1,Y_{0},X}(\delta |y_{0},x)& =\int_{\mathbb{R}
^{L}}h_{{\rm TT}}(\theta,\gamma_{2-L} ,x){\rm UCDITE}(\delta,\theta,\gamma_{2-L} ,x)
d\mathbb{P}_{Y_0,\overline{\Theta},\overline{\Gamma}_{2-L}|X}(y_0,\theta,\gamma_{2-L}|x)\ .
\end{align*}

\section*{Appendix 1: Modifications When Some of the Instruments are Discrete}
Because of its importance in applications, we consider the case where some instruments are discrete. 
For simplicity, we only detail the case with one additional instrument $B$ which is binary and use (A\ref{A5}).  We replace \eqref{m2}
by
\begin{equation}
D=\indic\left\{ \Theta -\alpha B -\Gamma_1(Z_1+g(Z_{2-L},B))-\sum_{l=2}^{L}\Gamma_lf_l(Z_l,B)>0\right\}\label{m2binary}
\end{equation}
where $\overline{\alpha}:=\alpha/\Gamma_1$, $g(Z_{2-L},B)$ and $f_l(Z_l,B)$ are integrable. Replace (A\ref{A2}), (A2'), (A\ref{A5}) by 
\begin{center}
\begin{tabular}{rl}
(AB2)\hspace{0.5cm} & $(Z,B)\ \perp\ (Y_d,\overline{\Theta},\overline{\Gamma}_{2-L})\ |X$ for $d=0,1$;\\
(AB4)\hspace{0.5cm} & $\forall(x,z_{2-L},b)\in \mathrm{supp}\left(\mathbb{P}_{X,Z_{2-L},B}\right)$, \\
&$
\mathrm{supp}\left( \mathbb{P}_{Z_1|Z_{2-L},B,X}(\cdot |z_{2-L},b,x)\right) \supseteq
\mathrm{supp}\left( \mathbb{P}_{\overline{\Theta}-\overline{\alpha} b-g(z_{2-L},b)-\sum_{l=2}^L\overline{\Gamma}_lf_l(z_l,b)|X}(\cdot |x)\right)$;
\end{tabular}
\end{center}
(N\ref{N0})-(N\ref{N3}) by (NB1)-(NB4) where we work with
$f_l(Z_l,B)$ for $l=2,\dots,L$, $f_l(Z_l)$ and $g(Z_{2-L},B)$; and the conditioning on $X$ by conditioning on $X$ and $B$.

Based on (AB2), (AB4), (NB1)-(NB4), the functions $g$ and $f_l$ for $l=2,\dots,L$ are identified and we can assume that the first-stage equation is 
\begin{equation}
D=\indic\left\{-\overline{\Theta} + \overline{\alpha} B +\sum_{l=2}^{L}\overline{\Gamma}_l\widetilde{Z}_L<\widetilde{Z}_1\right\}\label{m2binaryb}
\end{equation}
where $\widetilde{Z}_1:=-Z_1+g(Z_{2-L},B)$ and $\widetilde{Z}_l:=f_l(Z_l,B)$ for $l=2,\dots,L$.
We make Assumptions \ref{ass:I} and \ref{ass:Ib} on ${ \rm supp}(\mathbb{P}_{Z_{2-L}|B,X})$ where we also condition on $B$.

Under such assumptions, Theorem \ref{thrm:1} holds with $\mathbb{P}_{Y_d,\overline{\Theta}+\overline{\alpha}b,\overline{\Gamma}_{2-L},X}$ for $d=0,1$ and $b=0,1$. Thus, UCATE and UCQTE  are identified when we replace conditioning on $\overline{\Theta},\overline{\Gamma}_{2-L},X$ by $\overline{\Theta}+\overline{\alpha}b,\overline{\Gamma}_{2-L},X$ for $(\theta,\gamma_{2-L},b,x)$ in $\mathrm{supp}\left(\mathbb{P}_{\overline{\Theta}+\overline{\alpha}B,B,\overline{\Gamma}_{2-L}|B,X}\right)$. This gives, for $b=0,1$, two formulas for the effects which are weighted integrals.

UCDITE is identified when we replace (AB2) by
$(Z,B)\ \perp (Y_0,Y_1,\overline{\Theta},\overline{\Gamma}_{2-L})\ |X$, Assumption \ref{ass:condindep2} by
$Y_0 \perp Y_1\ |\overline{\Theta}+\overline{\alpha}b,\overline{\Gamma}_{2-L},X$, and Assumption \ref{ass:condindep3} \eqref{assi} by
$Y_0\perp\Delta\ |\overline{\Theta}+\overline{\alpha}b,\overline{\Gamma}_{2-L},X$. Without appealing to such assumptions, we obtain intersection bounds, for example, 
for $x$ in $\mathrm{supp}(\mathbb{P}_{X})$ and $(y_{0},y_{1})\in \mathbb{R}^{2}$,
\begin{equation*}
F_{Y_0,Y_1|X}^{L}(y_{0},y_{1}|x)\leq F_{Y_{0},Y_{1}|X}(y_{0},y_{1}|x)\leq
F_{Y_0,Y_1|X}^{U}(y_{0},y_{1}|x)
\end{equation*}
where, defining $F_{Y_0,Y_1|\overline{\Theta}+\overline{\alpha}b,\overline{\Gamma}_{2-L},X}^{L}(y_{0},y_{1}|\theta,\gamma_{2-L},x)$ and $F_{Y_0,Y_1|\overline{\Theta}+\overline{\alpha}b,\overline{\Gamma}_{2-L},X}^{U}(y_{0},y_{1}|\theta,\gamma_{2-L},x)$ like $F_{Y_0,Y_1|\overline{\Theta},\overline{\Gamma}_{2-L},X}^{L}(y_{0},y_{1}|\theta,\gamma_{2-L},x)$ and $F_{Y_0,Y_1|\overline{\Theta},\overline{\Gamma}_{2-L},X}^{U}(y_{0},y_{1}|\theta,\gamma_{2-L},x)$, 
\begin{align*}
F_{Y_0,Y_1}^{L}(y_{0},y_{1}|x)
& :=\max_{b=0,1}\int_{\mathbb{R}^{L}}F_{Y_0,Y_1|\overline{\Theta}+\overline{\alpha}b,\overline{\Gamma}_{2-L},X}^{L}(y_{0},y_{1}|\theta,\gamma_{2-L},x)
d\mathbb{P}_{\overline{\Theta}+\overline{\alpha}b,\overline{\Gamma}_{2-L}|X}(\theta,\gamma_{2-L}|x)\\
F_{Y_0,Y_1}^{U}(y_{0},y_{1}|x)
& :=\min_{b=0,1}\int_{\mathbb{R}^{U}}F_{Y_0,Y_1|\overline{\Theta}+\overline{\alpha}b,\overline{\Gamma}_{2-L},X}^{U}(y_{0},y_{1}|\theta,\gamma_{2-L},x)
d\mathbb{P}_{\overline{\Theta}+\overline{\alpha}b,\overline{\Gamma}_{2-L}|X}(\theta,\gamma_{2-L}|x)\ .
\end{align*}

Alternatively, effects conditional on $\{\overline{\Theta}=\theta,\overline{\alpha}=a,\overline{\Gamma}_{2-L}=\gamma_{2-L},X=x\}$ 
are identified under Assumption \ref{ass:condindep3}, replacing $Y_0,\Delta$ by $\overline{\Theta},\overline{\alpha}$ and without conditioning on $\overline{\Theta}$.

\section*{Appendix 2: Variance of Treatment Effects}
We now assume that $\mathbb{E}\left[Y_{0}^2+Y_{1}^2 \right] <\infty$, take $(\theta,\gamma_{2-L},x)\in\mathrm{supp}(\mathbb{P}_{\overline{\Theta},\overline{\Gamma}_{2-L},X})$, and use 
$$
V_1^U(\theta,\gamma_{2-L},x):=2\left(\mathbb{E}\left[\left.Y_0^2\right|\overline{\Theta}=\theta,\overline{\Gamma}_{2-L}=\gamma_{2-L},X=x\right]+\mathbb{E}\left[\left.Y_1^2\right|\overline{\Theta}=\theta,\overline{\Gamma}_{2-L}=\gamma_{2-L},X=x\right]\right)$$
\begin{align*}
V_2^U(\theta,\gamma_{2-L},x)&:=
\mathbb{E}\left[\left.Y_0^2\right|\overline{\Theta}=\theta,\overline{\Gamma}_{2-L}=\gamma_{2-L},X=x\right]+\mathbb{E}\left[\left.Y_1^2\right|\overline{\Theta}=\theta,\overline{\Gamma}_{2-L}=\gamma_{2-L},X=x\right]\\
&\quad-2\mathbb{E}\left[\left.Y_0\right|\overline{\Theta}=\theta,\overline{\Gamma}_{2-L}=\gamma_{2-L},X=x\right]\mathbb{E}\left[\left.Y_1\right|\overline{\Theta}=\theta,\overline{\Gamma}_{2-L}=\gamma_{2-L},X=x\right]\ .
\end{align*}
Under the assumptions of Theorem \ref{thrm:1} and assuming 
$$
 \mathbb{E}[(Y_0-\mathbb{E}[Y_0|\overline{\Theta},\overline{\Gamma}_{2-L},X])(Y_1-\mathbb{E}[Y_1|\overline{\Theta},\overline{\Gamma}_{2-L},X])|\overline{\Theta},\overline{\Gamma}_{2-L},X]\ge0\ a.s.
$$
for \eqref{ifposcorr} below, we easily obtain that
\begin{align}
Var\left(\left.\Delta\right|\overline{\Theta}=\theta,\overline{\Gamma}_{2-L}=\gamma_{2-L},X=x\right)&\leq
V_1^U(\theta,\gamma_{2-L},x)\notag\\
Var\left(\left.\Delta\right|\overline{\Theta}=\theta,\overline{\Gamma}_{2-L}=\gamma_{2-L},X=x\right)&\leq
V_2^U(\theta,\gamma_{2-L},x)\label{ifposcorr}
\end{align}

These upper bounds yield unconditional bounds when taking  expectation of the above inequalities evaluated at the random vectors. 

Based on the fact conditional variance of the sum of two uncorrelated variables is the sum of the conditional variances,
replacing Assumption \ref{ass:condindep3} by the weaker assumption
\begin{equation*}
 \mathbb{E}[(Y_0-\mathbb{E}[Y_0|\overline{\Theta},\overline{\Gamma}_{2-L},X])(\Delta-\mathbb{E}[\Delta|\overline{\Theta},\overline{\Gamma}_{2-L},X])|\overline{\Theta},\overline{\Gamma}_{2-L},X]=0\ a.s.
 \end{equation*}
yields
$$Var(\Delta|\overline{\Theta},\overline{\Gamma}_{2-L},X)=Var(Y_1|\overline{\Theta},\overline{\Gamma}_{2-L},X)-Var(Y_0|\overline{\Theta},\overline{\Gamma}_{2-L},X)\ .$$

\section*{Appendix 3: Proofs}
\noindent {\bf Proof of Theorem \ref{thrm:functions}.}
We have
$\mathbb{E}[W|X]=-\mathbb{E}[\overline{\Theta}|X]$,
while for $l=2,\dots,L$,
$\mathbb{E}[W|Z_l,X]=-\mathbb{E}[\overline{\Theta}|X]+\mathbb{E}[\overline{\Gamma}_l|X]f_l(Z_l)$,
thus
$f_l(Z_l)=\mathbb{E}[W|Z_l,X]-\mathbb{E}[W|X]$.
This identifies the functions $f_l$ on $\mathrm{supp}(\mathbb{P}_{Z_l})$ for $l=2,\dots,L$. 
Finally, $g$ is identified using
$$g(Z_{2-L})=\mathbb{E}[W|Z_{2-L},X]-\mathbb{E}[W|X]-\sum_{l=2}^Lf_l(Z_l)\ .$$

\noindent {\bf Proof of Theorem \ref{thrm:1}.} (A\ref{A2}) implies that for $(z,x)\in\mathrm{supp}(\mathbb{P}_{Z,X})$, $y_1\in\mathbb{R}$
\begin{align}
\mathbb{E}\left[\left.\indic\{Y\le y_1\}D\right|\overline{Z}=z,X=x\right]&= \mathbb{E}\left[\left.\indic\{Y_1\le y_1\}D\right|\overline{Z}=z,X=x\right]\notag\\
&=\mathbb{E}\left[\left.\indic\{Y_1\le y_1\}\indic\left\{ -\overline{\Theta}+\sum_{l=2}^{L}\overline{\Gamma}_lz_l<z_1\right\}\right|X=x\right]\notag\ .
\end{align}
Using (A\ref{A3}) and (A\ref{A5}), this yields that $F_{Y_1,-\overline{\Theta}+\sum_{l=2}^{L}\overline{\Gamma}_lz_l|X}(\cdot|x)$ is identified. Thus, for every $t,s$ in $\mathbb{R}$, $x\in \mathrm{supp}(\mathbb{P}_{X})$, and $z_{2-L}$ in $\mathrm{supp}(\mathbb{P}_{Z_{2-L}|X}(\cdot|x))$,
$\mathbb{E}\left[\left.\exp\left(itY_1-is\overline{\Theta}+is\sum_{l=2}^{L}\overline{\Gamma}_lz_l\right)\right|X=x\right]$ is identified.  Under Assumption \ref{ass:Ib}, (A\ref{S1}) implies that we can extend as an entire function
$z_{2-L}\to \mathbb{E}\left[\left.\exp\left(itY_1-is\overline{\Theta}+is\sum_{l=2}^{L}\overline{\Gamma}_lz_l\right)\right|X=x\right]$ and (S\ref{S2}) that for every $t,s\in\mathbb{R}$ and $z_{2-L}\in\mathbb{R}$,  $\mathbb{E}\left[\left.\exp\left(itY_1-is\overline{\Theta}+is\sum_{l=2}^{L}\overline{\Gamma}_lz_l\right)\right|X=x\right]$ is identified.  Thus,  $(t,s,z_{2-L})\in \mathbb{R}^2\times (\mathbb{R}\setminus\{0\})^{L-1}\to \mathbb{E}\left[\left.\exp\left(itY_1+is\overline{\Theta}+i\sum_{l=2}^{L}\overline{\Gamma}_lz_l\right)\right|X=x\right]$ is identified. 
We obtain directly this result if we make Assumption \ref{ass:I} instead of Assumption \ref{ass:Ib}. Now, by continuity, the 
Fourier transform of $\mathbb{P}_{Y_1,\overline{\Theta},\overline{\Gamma}_{2-L}|X}(\cdot|x)$ is identified everywhere. The injectivity of the Fourier transform of measures allows to conclude.

If we maintain (A4') instead of (A\ref{A5}), the above argument yields that $\mathbb{P}_{Y_d,-\overline{\Theta}+\sum_{l=2}^{L}\overline{\Gamma}_lz_l}$ for $j=0,1$ is identified on $\mathbb{R}\times\mathrm{supp}(\mathbb{P}_{Z_1|Z_{2-L},X}(\cdot|z_{2-L},x))$ for every $(z_{2-L},x)$ in $\mathrm{supp}(\mathbb{P}_{Z_{2-L},X})$. (A4'\ref{SLS1}) and (A4'\ref{SLS2}) now implies that, for every $t,s$ in $\mathbb{R}$,  $(z_{2-L},x)$ in $\mathrm{supp}\left(\mathbb{P}_{Z_{2-L},X}\right)$ and $z_1$ in $\mathrm{supp}\left(\mathbb{P}_{Z_1|\overline{Z}_{2-L},X}(\cdot|z_{2-L},x)\right)$, $\mathbb{E}\left[\left.\exp\left(itY_1-is\overline{\Theta}+iz_1\sum_{l=2}^{L}\overline{\Gamma}_lz_l\right)\right|X=x\right]$ is identified. Now, we can either conclude under Assumption \ref{ass:I} or Assumption \ref{ass:Ib}.

Identification of $F_{Y_0,\overline{\Theta},\overline{\Gamma}_{2-L}}$ is obtained in the same way replacing $D$ by $1-D$.\vspace{0.3cm}

\noindent {\bf Proof of Theorem \ref{prop:ucate}}. This follows from the fact that, based on \eqref{eq:info1}-\eqref{eq:info3}, we have
\begin{align*}
\mathbb{E}[Y_1-Y_0|\overline{\Gamma}_{2-L},\eta_1-\eta_0,X]
&=\mathbb{E}\left[\left.\mathbb{E}[Y_1-Y_0|\mathcal{I}]\right|
\overline{\Gamma}_{2-L},\eta_1-\eta_0,X\right]\\
&=\mathbb{E}\left[\left.\mathbb{E}[Y_1-Y_0|\overline{\Gamma}_{2-L},X]+\eta_1-\eta_0\right|
\overline{\Gamma}_{2-L},\eta_1-\eta_0,X\right]\\
&=\mathbb{E}[Y_1-Y_0|\overline{\Gamma}_{2-L},X]+\eta_1-\eta_0=\mathbb{E}[\Delta|\mathcal{I}]\ .
\end{align*}
Now, because $\mathbb{E}[\Delta|\mathcal{I}]$ is not $\sigma(\overline{\Gamma}_{2-L},X)$-measurable and we assume (U\ref{U2}), for every $(\gamma_{2-L},x)\in \mathrm{supp}(\mathbb{P}_X)$, there are two draws $\theta_1\ne \theta_2$ from $F_{\overline{\Theta}|\overline{\Gamma}_{2-L},X}(\cdot|\gamma_{2-L},x)$ such that $$\Gamma_1(\gamma_{2-L},x)(\theta_1-\theta_2)=\mathrm{UCATE}(\theta_1,\gamma_{2-L},x)-\mathrm{UCATE}(\theta_2,\gamma_{2-L},x)\ne0\ .$$
Thus $\Gamma_1(\gamma_{2-L},x)$ and $\Gamma_0(\gamma_{2-L},x)=\mathrm{UCATE}(\theta_1,\gamma_{2-L},x)-\Gamma_1(\gamma_{2-L},x)\theta_1$
are identified. Stated otherwise, conditioning on a fixed value of $(\overline{\Gamma}_{2-L},X)$, $\Gamma_0$ and $\Gamma_1$ are constant and 
parameters in a regression without error. Thus EAR is uniquely determined once we condition on a fixed value of $(\overline{\Theta},\overline{\Gamma}_{2-L},X)$ and $Z$. 
This yields the result by integration.

\noindent {\bf Proof of Theorem \ref{theorem:BoundsMU}.} The first inequality is based on the Frechet-Hoeffding bounds applied to the conditional CDFs which are CDFs as well. The second one is a consequence of the Makarov bounds (Makarov (1981)) for the same reason. Sharpness is discussed in Firpo and Ridder (2008). The last two inequalities are obtained by simple manipulations.\vspace{0.3cm}

\noindent {\bf Proof of Corollary \ref{cor:1}.} 
The bounds are obtained by taking expectations of the inequalities in Theorem \ref{theorem:BoundsMU} evaluated at the random vectors.\vspace{0.3cm}  

\noindent {\bf Proof of Theorem \ref{thrm:CondIndep}.} 
Assumption \ref{ass:condindep1} implies that
$$F_{Y_0,Y_1|\overline{\Gamma}_{2-L},X}(y_1,y_0|\gamma_{2-L} ,x)=F_{Y_0|\overline{\Gamma}_{2-L},X}(y_0|\gamma_{2-L} ,x)F_{Y_1|\overline{\Gamma}_{2-L},X}(y_1|\gamma_{2-L} ,x)\ ,$$ 
similarly, Assumption \ref{ass:condindep2} implies that
$$F_{Y_0,Y_1|\overline{\Theta},\overline{\Gamma}_{2-L},X}(y_0,y_1|\theta,\gamma_{2-L},x)=F_{Y_0|\overline{\Theta},\overline{\Gamma}_{2-L},X}(y_0|\theta,\gamma_{2-L},x)F_{Y_1|\overline{\Theta},\overline{\Gamma}_{2-L},X}(y_1|\theta,\gamma_{2-L},x)$$ and the right-hand sides are identified.

Let us now consider Assumption \ref{ass:condindep3}. Based on Theorem \ref{thrm:1}, $F_{Y_0|\overline{\Theta},\overline{\Gamma}_{2-L},X}(\cdot|\theta,\gamma_{2-L},x)$ and $F_{Y_0+\Delta | \overline{\Theta},\overline{\Gamma}_{2-L},X}(\cdot|\theta,\gamma_{2-L},x)$ are identified. Estimating $F_{\Delta | \overline{\Theta},\overline{\Gamma}_{2-L},X}(\cdot|\theta,\gamma_{2-L},x)$ and thus $F_{Y_0,Y_1| \overline{\Theta},\overline{\Gamma}_{2-L},X}(\cdot|\theta,\gamma_{2-L},x)$ is a deconvolution problem under \eqref{assi}.\vspace{0.3cm}

\noindent {\bf Proof of Theorem \ref{thrm:ex}.} 
This is a consequence of theorems \ref{thrm:CondIndep} (and its proof) and \ref{prop:ucate}.


\begin{thebibliography}{99}
\bibitem{AAI} \textsc{Abadie, A., J. Angrist, and G. Imbens} (2002):
{``Instrumental Variables Estimates of the Effect of Subsidized
Training on the Quantiles of Trainee Earnings"}.
\emph{Econometrica}, \textbf{70}, 91--117.

\bibitem{AH} \textsc{Abbring, J. H., and J. J. Heckman} (2007):
{``Econometric Evaluation of Social Programs, Part III:
Distributional Treatment Effects, Dynamic Treatment Effects, Dynamic
Discrete Choice, and General Equilibrium Policy Evaluation"}.
\emph{Handbook of Econometrics}, J.J. Heckman and E.E. Leamer
(eds.), Vol. 6, North Holland, Chapter 72.

\bibitem{BH} \textsc{Beran, R., and P. W. Millar} (1994):
{``Minimum Distance Estimation in Random Coefficients Regression Models"}. \emph{Annals of Statistics}, \textbf{22}, 1976--1992.

\bibitem{BFH} \textsc{Beran, R., A. Feuerverger, and P. Hall} (1996):
{``On Nonparametric Estimation of Intercept and Slope in Random
Coefficients Regression"}. \emph{Annals of Statistics}, \textbf{24},
2569--2592.

\bibitem{BM} \textsc{Bj\"{o}rklund, A., and R. Moffitt} (1987): {
\textquotedblleft The Estimation of Wage and Welfare Gains in Self-Selection
Models"}, \emph{Review of Economics and Statistics}, \textbf{69}, 42--49.

\bibitem{CHH} \textsc{Carneiro, P., K. T. Hansen, and J. Heckman} (2003):
{``Estimating Distributions of Treatment Effects With an Application
to the Return to Schooling and Measurement of the Effect of
Uncertainty on College Choice"}. \emph{International Economic
Review}, \textbf{44}, 361--422.

\bibitem{CHH} \textsc{Carneiro, P., J. Heckman, and E. Vytlacil} (2010): {
\textquotedblleft Evaluating Marginal Policy Changes and the Average Effect
of Treatment for Individuals at the Margin"}. \emph{Econometrica}, \textbf{78},
377--394.

\bibitem{CH1} \textsc{Chernozhukov, V., and C. Hansen} (2005): {``An IV Model of Quantile Treatment Effects"}. \emph{\ Econometrica}, \textbf{73}, 245--261.

\bibitem{CunhaHeckman}
\textsc{Cunha, F., and J. Heckman} (2007): {``Identifying and Estimating the Distributions of Ex Post and Ex Ante Returns to Schooling"}. \emph{Labour Economics}, \textbf{14}, 870--893.

\bibitem{Devroye}
\textsc{Devroye, L.} (1989): {``Consistent Deconvolution in Density Estimation"}. \emph{The Canadian Journal of Statistics}, 
\textbf{17}, 235--239.

\bibitem{FZ} \textsc{Fan, Y., E. Guerre and D. Zhu} (2014): {``Partial Identification and Confidence Sets for Functional of the Joint Distribution of Potential Outcomes"}. Working paper.

\bibitem{FP} \textsc{Fan, Y., and S. S. Park} (2010): {\textquotedblleft
Sharp Bounds on the Distribution of Treatment Effects and Their Statistical
Inference"}. \emph{Econometric Theory}, \textbf{26}, 931--951.

\bibitem{FR} \textsc{Firpo, S., and G. Ridder} (2008): {``Bounds on
Functionals of the Distribution of Treatment Effects"}. Working paper.

\bibitem{FG} \textsc{Fox, J., and A. Gandhi} (2011): {\textquotedblleft A
Simple Nonparametric Approach to Estimating the Distribution of Random
Coefficients in Structural Models"}. Working Paper.

\bibitem{GG} \textsc{Gaillac, C., and E. Gautier} (2015):
{\textquotedblleft Estimation of the Distribution of Random Coefficients with Bounded Regressors"}. Working paper.

\bibitem{GK} \textsc{Gautier, E., and Y. Kitamura} (2013):
{\textquotedblleft Nonparametric Estimation in Random Coefficients Binary
Choice Models"}. \emph{Econometrica}, \textbf{81}, 581--607.

\bibitem{HSC} \textsc{Heckman, J. J., J. Smith, and N. Clements} (1997): {\textquotedblleft Making The Most Out Of Programme Evaluations and Social Experiments: Accounting For Heterogeneity in Programme Impacts"}. \emph{Review of Economic Studies}, \textbf{64}, 487--635.

\bibitem{HV} \textsc{Heckman, J. J., and E. Vytlacil} (2001) {
\textquotedblleft Policy Relevant Treatment Effects"}. \emph{American
Economic Review Papers and Proceedings}, \textbf{91}, 107--111.

\bibitem{HV2} \textsc{Heckman, J. J., and E. Vytlacil} (2005): {
\textquotedblleft Structural Equations, Treatment Effects, and Econometric
Policy Evaluation"}. \emph{Econometrica}, \textbf{73}, 669--738.

\bibitem{HKM} \textsc{Hoderlein, S., J. Klemel\"{a}, and E. Mammen} (2010): {
\textquotedblleft Analyzing the Random Coefficient Model Nonparametrically"}. \emph{Econometric Theory}, \textbf{26}, 804--837.
 
\bibitem{IT} \textsc{Ichimura, H., and T. S. Thompson} (1998): {
\textquotedblleft Maximum Likelihood Estimation of a Binary Choice Model
with Random Coefficients of Unknown Distribution"}. \emph{Journal of
Econometrics}, \textbf{86}, 269--295.

\bibitem{IA} \textsc{Imbens, G. W., and J. D. Angrist} (1994): {
\textquotedblleft Identification and Estimation of Local Average Treatment
Effects"}. \emph{Econometrica}, \textbf{62}, 467--475.

\bibitem{IN} \textsc{Imbens, G. W., and W. K. Newey} (2009): {
\textquotedblleft Identification and Estimation of Triangular Simultaneous
Equations Models Without Additivity Corresponding"}. \emph{Econometrica},
\textbf{77}, 1481--1512.

\bibitem{Kasy}
\textsc{Kasy, M.} (2011): {\textquotedblleft Identification in Triangular Systems Using Control Functions"}. \emph{Econometric Theory}, \textbf{27}, 663--671.

\bibitem{Kitagawa}
\textsc{Kitagawa, T.} (2009): {\textquotedblleft
Identification Region of the Potential Outcome Distributions under Instrument Independence"}. Working paper.

\bibitem{Klein} \textsc{Klein, T.} (2010): {\textquotedblleft Heterogeneous
Treatment Effects: Instrumental Variables Without Monotonicity?"}. \emph{
Journal of Econometrics}, \textbf{155}, 99--116.

\bibitem{Lewbel2} \textsc{Lewbel, A.} (2007): {\textquotedblleft Endogenous
Selection or Treatment Model Estimation"}. \emph{Journal of Econometrics},
\textbf{141}, 777--806.

\bibitem{LP}
\textsc{Lewbel, A. and K. Pendakur} (2012):
{\textquotedblleft Generalized Random Coefficients With Equivalence Scale Applications"}. Working paper.

\bibitem{Makarov} \textsc{Makarov, G. D.} (1981): {``Estimates of the Distribution Function of a Sum of Two Random Variables when the Marginal Distributions are Fixed"}. \emph{Theory of Probability and its Applicatons}, \textbf{26}, 803--806.

\bibitem{NeweyPowellVella}
\textsc{Newey, W. K., Powell, J. L. and F. Vella} (1999):
{\textquotedblleft Nonparametric Estimation of Triangular
Simultaneous Equations Models"}. \emph{Econometrica}, \textbf{67}, 565--603.

\bibitem{ST}
\textsc{Small and Tan} (2007):
{\textquotedblleft A Stochastic Monotonicity Assumption for the Instrumental Variables Method"}. Working paper.

\bibitem{Vytlacil}
\textsc{Vytlacil, E.} (2002): {\textquotedblleft Independence, Monotonicity, 
and Latent Index Models: An Equivalence Result"}. \emph{Econometrica}, \textbf{70}, 331--341.

\end{thebibliography}
\end{document}